\documentclass[12pt]{amsart}
\usepackage{graphicx}
\usepackage{amssymb,amstext,amsfonts,amscd}

\newtheorem{theorem}{Theorem}
\newtheorem{lemma}[theorem]{Lemma}

\numberwithin{equation}{section}

\newcommand{\bn}{\mathbb{N}}
\newcommand{\bz}{\mathbb{Z}}
\newcommand{\br}{\mathbb{R}}
\newcommand{\bc}{\mathbb{C}}
\newcommand{\bp}{\mathbb{P}}
\newcommand{\bq}{\mathbb{Q}}
\newcommand{\bff}{\mathbb{F}}
\newcommand{\M}{\mathcal{M}}

\newcommand{\bhh}{\mathbb{H}}

\newcommand{\rk}{\mathop\mathrm{rk}}

\newcommand{\Proof}{\noindent{\it{Proof.}}\ \ }
\newcommand{\QED}{\ \ $\Box$}
\newcommand{\wP}{\widetilde{P}}

\title{Remarks on connected components of moduli of real polarized K3 surfaces}

\author{Viacheslav V. Nikulin}

\address{Department of Pure Mathematics, The University of Liverpool,
Liverpool, L69 3BX, United Kingdom\\
Steklov Mathematical Institute, ul. Gubkina 8, Moscow, GSP-1,
Russia} \email{vnikulin@liv.ac.uk, vvnikulin@list.ru}


\begin{document}
\pagestyle{plain}


\begin{abstract}
We have finalized the old (1979) results from \cite{Nikulin79}
about enumeration of connected components of moduli or real
polarized K3 surfaces.

As an application, using recent results of \cite{NikulinSaito05}
(see also math.AG/0312396), we completely classify real
polarized K3 surfaces which are deformations of real hyper-elliptically
polarized K3 surfaces. This could be important in some questions,
because real hyper-elliptically polarized K3 surfaces can be 
constructed explicitly. 
\end{abstract}

\maketitle

\tableofcontents

\section{Introduction}

In \cite{Nikulin79}, enumeration of connected components of moduli of
real polarized K3 surfaces $(X,P^\prime)$ had been considered. We review
these results in Section \ref{old}.

Using Global Torelli Theorem for K3 surfaces due to 
Piatetski\u i-Shapiro and Shafarevich 
\cite{PS71} and epimorphicity of Torelli map for K3 surfaces 
due to Vic. Kulikov 
\cite{Kulikov77}, it had been shown in \cite{Nikulin79} that the connected
component of moduli is determined by the isomorphism class of the action
of the anti-holomorphic involution
$\varphi$ in $H_2(X(\bc),\bz)$ with considering the corresponding polarization
$P^\prime\in H_2(X(\bc),\bz)$ which satisfies $\varphi (P^\prime)=-P^\prime$.
This reduces enumeration of connected components to a purely arithmetic
problem.

To solve this problem, in \cite{PS71} the invariants
\begin{equation}
(r,a,\delta_\varphi;k,n,\delta_P,\delta_{\varphi P}) \label{IntgenK3}
\end{equation}
of the genus of this action
had been introduced and classified. They have very clear
geometric meaning. For example, $k\in \bn$ is defined by the condition
that $P=P^\prime/k$ is primitive in the Picard lattice, and $n=P^2$ is
the primitive degree, $n>0$ and $n\equiv 0\mod 2$. The invariants
$(r,a,\delta_\varphi)$ determine
the topological type of $X(\br)$ which is an orientable compact surface. The
invariant $\delta_{\varphi P}\in \{0,1\}$, and $\delta_{\varphi P}=0$ if
and only if $X(\br)\sim P\mod 2$ in $H_2(X(\bc),\bz)$.
All possible invariants $(r,a,\delta_\varphi)$ are given in Figure
\ref{Sgraph} and describe connected components of moduli of real K\"ahlerian
(i. e. without polarization) K3 surfaces. See Sect. \ref{new} for
other details.

The invariants \eqref{IntgenK3} give main invariants of connected components
of moduli of real polarized K3 surfaces. It had been shown in
\cite{Nikulin79}, that the invariants \eqref{IntgenK3}
``almost always'' determine the connected component of moduli uniquely, i. e.
the moduli space of real polarized K3 surfaces with these invariants
is connected.
More exactly, by \cite{Nikulin79}, this is true if $r\le 18-\delta_P$.
Here the invariant $r\in \bn$ takes values in
the range $1\le r\le 20$ and $\delta_P\in \{0,1\}$. Thus, only for very few
cases when $r=18$ or $r=19,\ 20$ (look at Figure \ref{Sgraph})
the connectedness was not known after the results of \cite{Nikulin79}.

In Sect. \ref{new} we finalize this uniqueness result.
We show that the same uniqueness is valid if $r\le 18$.
To prove this, we use results of Miranda and Morrison
\cite{MirandaMorrison85}, \cite{MirandaMorrison86} and of D. James
\cite{James96} where the general analogue of Witt's  
theorem had been proved for
indefinite even lattices of the rank at least three. 
In \cite{Nikulin79} its particular 
case (Theorem 1.14.2 in \cite{Nikulin79}) had been  used.
Real K3 surfaces which satisfy
the condition $r\le 18$ can be described purely topologically:
$$
X(\br)\not=T_1\amalg (T_0)^8,\ (T_0)^9,\ (T_0)^{10}.
$$
Here $T_g$ denotes a compact orientable surface of the genus $g$. Thus,
if the topological type is different from these three types, then
the connected component of moduli of real polarized K3 surfaces is
determined uniquely by the genus invariants \eqref{IntgenK3}.

If $X(\br)=T_1\amalg (T_0)^8$, $(T_0)^9$ or $(T_0)^{10}$, then
the uniqueness result above surely is not valid in general,
for any primitive
degree $n$, and one has to introduce some additional to \eqref{IntgenK3}
invariants of connected components of moduli.
We do it in Theorems \ref{theorem1911} --- \ref{theorem202} introducing
the additional invariants 
which determine the connected component of moduli uniquely. Here we use
the discriminant forms technique which had been developed in \cite{Nikulin79}.

For real polarized K3 surfaces having the topological type
$X(\br)=T_1\amalg (T_0)^8$ or $(T_0)^9$ some connected components of moduli
(and their K3) are especially remarkable and are called
{\it standard} or {\it different from standard only over $2$.} They are
determined by the especially simple these additional invariants, and
are important for our further considerations. See their
definition after formulations of the theorems
\ref{theorem1911}---\ref{theorem202}.

In Sect. \ref{sechyppolK3}, using recent results of \cite{NikulinSaito05},
we apply these improvements of the old results of \cite{Nikulin79} to answer
the following interesting question:
{\it Which real polarized K3 surfaces are deformations of
real hyper-elliptically polarized K3 surfaces $(X,\wP)$?} Here $(X,\wP)$
is hyper-elliptically polarized when the
general curve of the corresponding complete linear
system $|\wP|$ is hyper-elliptic, and then the linear system gives
a double covering. Roughly speaking, when a connected component of moduli
of real polarized K3 surfaces contains a hyper-elliptically
polarized K3 surface?

This could be important in some questions because real hyper-elliptically
polarized K3 surfaces can be constructed explicitly
as double coverings of the relatively minimal rational surfaces $Y=\bp^2$
or $Y=\bff_n$ for $n=0,\,1,\,2,\,4$ ramified
in a non-singular curve $A\in |-2K_Y|$ where $K_Y$ denotes the canonical class
of $Y$. The corresponding complete linear system $|\wP|$ is then
the preimage of some standard complete linear system of $Y$.

In \cite{NikulinSaito05} all possible genus invariants \eqref{IntgenK3}
of these deformations had been described (in \cite{NikulinSaito05} they
were called as deformations of general K3 double rational scrolls which is
equivalent). Using these results and the results above about
enumeration of connected component of moduli of real polarized K3 surfaces,
we can classify these K3 surfaces completely.

We get the following result where the first statement (i)
had been obtained in \cite{NikulinSaito05}.

\begin{theorem}
\label{InttheoremhypK3}
A real polarized K3 surface $(X,P^\prime)$ is
a deformation of a general real K3 double rational scroll
(equivalently, of a real hyper-elliptically polarized K3 surface) if and
only if one of conditions (i)---(iv) below satisfies:

(i) The primitive degree $n=2$ or $4$ (see \cite{NikulinSaito05}).

(ii) The primitive degree $n\ge 6$, and
$X(\br)\not=T_1\amalg (T_0)^8,\ (T_0)^9,\ (T_0)^{10}$, and
$X(\br)\not\sim P\mod 2$ in $H_2(X(\bc),\bz)$ if $X(\br)=(T_0)^k$.

(iii) The primitive degree $n\ge 6$, and
$X(\br)=T_1\amalg (T_0)^8$, and $(X,P^\prime)$ is standard if
$n\equiv 0,\ 2\mod 8$, and $(X,P^\prime)$ is different from
a standard only over $2$ if $n\equiv 4,\ 6\mod 8$.

(iv) The primitive degree $n\ge 6$, and
$X(\br)=(T_0)^9$, and $X(\br)\not\sim P\mod 2$ in $H_2(X(\bc),\bz)$, and
$(X,P^\prime)$ is standard.
\end{theorem}

\label{introduction}

\section{Enumeration of connected components of moduli of real polarized K3
surfaces}\label{secmoduliK3}

Here we finalize results of \cite{Nikulin79} about description of
connected components of moduli of real polarized K3 surfaces.

\subsection{Reminding of known results about connected components
of moduli of real polarized K3 surfaces} \label{old}
Here we review results of \cite{Nikulin79}.

Let
$(X,P^\prime)$ be a real polarized K3 surface. Here $X$ is an
algebraic K3 surface, and $P^\prime$ a very ample divisor class on
$X$, defined over the field $\br$ of real numbers. I. e.
$(X,P^\prime)$ is a complex polarized K3 surface together with an
anti-holomorphic involution $\varphi$  of $X$ such that
$\varphi^\ast (P^\prime)=-P^\prime$. We want to describe connected
components of moduli of the pairs $(X,P^\prime)$.

Let $L=H_2(X(\bc),\bz)$ be the homology lattice of $X$ with the
intersection pairing. It is an even unimodular lattice of
signature $(3,19)$. This characterizes the lattice up to
isomorphisms. The polarization $P^\prime\in L$ is an element of
$L$ with the $(P^\prime)^2>0$. The anti-holomorphic involution
$\varphi$ acts in $L$, and $\varphi(P^\prime)=-P^\prime$. The
triplet $(L,\varphi,P^\prime)$ considered up to natural
isomorphisms is called {\it the polarized integral K3 involution
corresponding to the real polarized K3 surface $(X,P^\prime)$.}
Here another triplet $(\widetilde{L},\widetilde{\varphi},
\widetilde{P^\prime})$ is isomorphic to $(L,\varphi,P^\prime)$ if
there exists an isomorphism $f:L\to\widetilde{L}$ of lattices (i.
e. preserving the intersection pairing) such that
$f(P^\prime)=\widetilde{P^\prime}$ and
$\widetilde{\varphi}f=f\varphi$. Further we denote by $L^\varphi$
and $L_\varphi$ the eigenvalue $1$ and $-1$ parts respectively of
the action of $\varphi$ in $L$. The integral polarized K3
involution $(L,\varphi,P^\prime)$ satisfies the conditions:
$L^\varphi$ is hyperbolic (i. e. it has exactly one positive
square), and $P^\prime\in L_\varphi$ (e.g. see \cite{Kharlamov76},
\cite{Nikulin79}).

By Theorem 3.10.1 in \cite{Nikulin79}, we have the following main
result which is based on Global Torelli Theorem  \cite{PS71} and
epimorphicity of Torelli map \cite{Kulikov77} for K3 surfaces.

\begin{theorem} (see Theorem 3.10.1 in \cite{Nikulin79})
Connected components of moduli of real polarized K3 surfaces are
in one to one correspondence with isomorphism classes of integral
polarized involutions $(L,\varphi,P^\prime)$ such that $L$ is an
even unimodular lattice of signature $(3,19)$, $L^\varphi$ is
hyperbolic, $P^\prime\in L_\varphi$ (i. e.
$\varphi(P^\prime)=-P^\prime$) and $(P^\prime)^2>0$.
\label{intpolinvtheorem}
\end{theorem}

Thus, description of connected components of moduli of real
polarized K3 surfaces is equivalent to a purely arithmetic problem
of classification of the integral polarized involutions above.
Further we call them as integral polarized K3 involutions. To
solve this problem, in \cite{Nikulin79}, their invariants
\begin{equation}
(r,a,\delta_\varphi;k,n,\delta_P,\delta_{\varphi P}) \label{genK3}
\end{equation}
were introduced. Here $r=\rk L^\varphi\in \bn$;
$((L^\varphi)^\ast/L^\varphi)=(\bz/2\bz)^a$ where $a\ge 0$ is an
integer; $\delta_\varphi \in \{0,1\}$ is equal to $0$ if and only
if $x\cdot \varphi(x)\equiv 0\mod 2$ for any $x\in L$. They are
all invariants of the corresponding pair $(L,\varphi)$. Here $k\in
\bn$ is defined by the condition that $P=P^\prime/k$ is a
primitive element of $L$; here {\it the primitive degree}
$n=P^2=(P^\prime/k)^2$ is an even
natural number; here $\delta_{P}\in \{0,1\}$ is equal to $0$ if
and only if $P \cdot L_\varphi\equiv 0\mod 2$; here
$\delta_{\varphi P}\in \{0,1\}$ is equal to $0$ if and only if
$x\cdot \varphi(x)\equiv x\cdot P$ for any $x\in L$. The
invariants \eqref{genK3} give all {\it invariants of the genus} of
the corresponding integral polarized K3 involutions: for equal
invariants \eqref{genK3}, the corresponding integral polarized K3
involutions are isomorphic over $\br$ and the rings $\bz_p$ of
$p$-adic integers for any prime $p$.

See \cite{Kharlamov75a}, \cite{Kharlamov76} and \cite{Nikulin79}
about geometric meaning of the invariants \eqref{genK3}. We only
mention the following where we denote by $T_g$ an orientable
compact surface of the genus $g$. We have
\begin{equation}
X(\br)=
\begin{cases}
\emptyset &\text{if $(r,a,\delta_\varphi)=(10,10,0)$}\\
T_1\amalg T_1 &\text{if $(r,a,\delta_\varphi)=(10,8,0)$}\\
T_g\amalg (T_0)^k &\text{otherwise, where}\\
              &g=(22-r-a)/2,\ k=(r-a)/2
\end{cases};
\label{realcomponents}
\end{equation}
\begin{equation}
X(\br)\sim 0 \mod 2 \text{\ in\ } H_2(X(\bc);\bz) \label{realmod21}
\end{equation}
if and only if $\delta_{\varphi}=0$, and
\begin{equation}
X(\br)\sim P \mod 2 \text{\ in\ } H_2(X(\bc);\bz) \label{realmod22}
\end{equation}
if and only if $\delta_{\varphi P}=0$. Here $X(\br)=X(\bc)^\varphi$ is
the fixed points set for the corresponding anti-holomorphic
involution $\varphi$ on the complex K3 surface $X(\bc)$.

In (\cite{Nikulin79}, Theorem 3.4.3), the genus invariants
\eqref{genK3} of the integral polarized K3 involutions were
classified: One should set $l_{(+)}=3$, $l_{(-)}=19$, $t_{(+)}=1$
and $t_{(-)}=r-1$ in this theorem. We have

\begin{theorem}
 (see Theorem 3.4.3 in \cite{Nikulin79}) The
invariants \eqref{genK3} give complete genus invariants of integral
polarized K3 involutions.

There exists a real polarized K3 surface with the genus invariants
\eqref{genK3} if and only if the invariants satisfy the conditions
0.(1)--(7) and I. (1)---(19) listed below.

\label{geninvtheorem}
\end{theorem}

\noindent {\bf 0. Conditions on $(r,a,\delta_\varphi )$:}

\noindent (1) $1\le r\le 20$, $0\le a\le \min \{r,\ 22-r\}$;

\noindent (2) $r+a\equiv 0\mod 2$; if $\delta_\varphi=0$, then
$r\equiv 2\mod 4$;

\noindent (3) if $a=0$, then $\delta_\varphi=0$ and $r\equiv 2\mod
8$;

\noindent (4) if $a=1$, then $r\equiv 1,\,3\mod 8$;

\noindent (5) if $(a=2,\  r\equiv 6\mod 8)$, then
$\delta_\varphi=0$;

\noindent (6) if $(a=r,\ \delta_\varphi=0)$, then $r\equiv 2\mod
8$;

\noindent (7) if $(a=22-r, \delta_\varphi=0)$, then $r\equiv 2\mod
8$.

\medskip

\noindent{\bf I. Conditions on $n$, $\delta_P$, $\delta_{\varphi
 P}$:}

\noindent
{\bf General conditions:}

\noindent (1) $n>0$ and $n\equiv 0\mod 2$;

\noindent (2) if $(n\equiv 2\mod 4,\ \delta_P=0)$, then
$\delta_{\varphi}=1$;

\noindent (3) if $\delta_{\varphi P}=0$, then $(\delta_P=0,\
\delta_{\varphi }=1,\ r\equiv 2+n/2 \mod 4)$.

\medskip

\noindent {\bf Relations near the boundary $a=22-r$:}

\noindent (4) if $a=22-r$, then $\delta_P=0$;

\noindent (5) if $(a=22-r,\ \delta_{\varphi P}=0)$, then $r\equiv
2+n/2\mod 8$;

\noindent (6) if $(a=20-r,\ n\equiv 0\mod 4,\ \delta_P=1,\
\delta_\varphi=0)$, then $r\equiv 2\mod 8$.

\medskip

\noindent {\bf Relations near the boundary $a=0$:}

\noindent (7) if $a=0$, then $\delta_P=1$;

\noindent (8) if $(a=1,\ n\equiv 0\mod 4)$, then $\delta_P=1$;

\noindent (9) if $(a=1,\ \delta_P=0,\ n\equiv \pm 2\mod 8)$, then
$r\equiv 2\pm 1\mod 8$;

\noindent (10) if $(a=2,\ \delta_P=0,\ n\equiv \pm 2\mod 8)$, then
$r\equiv 2,\, 2\pm 2\mod 8$;

\noindent (11) if $(a=2,\ \delta_P=0,\ n\equiv 0\mod 8)$, then
$r\equiv 2 \mod 8$;

\noindent (12) if $(a=3,\ \delta_P=0,\ n\equiv 0\mod 8)$, then
$r\equiv 1,\,3\mod 8$;

\noindent (13) if $(a=2,\ \delta_P=0,\ n\equiv 4\mod 8,\ r\equiv
2\mod 8)$, then $\delta_{\varphi}=0$;

\noindent (14) if $(a=1,\ \delta_P=0)$, then $\delta_{\varphi
P}=0$;

\noindent (15) if $(a=2,\ \delta_P=0,\ n\equiv 4\mod 8,\ r\equiv
0\mod 4)$, then $\delta_{\varphi P}=0$;

\noindent (16) if $(a=3,\ \delta_P=0,\ n\equiv \pm 2\mod 8,\
r\equiv 2\pm 5\mod 8)$, then $\delta_{\varphi P}=0$;

\noindent (17) if $(a=2,\ \delta_P=0,\ n\equiv 0\mod 8,\ r\equiv
2\mod 8,\ \delta_{\varphi}=1)$, then $\delta_{\varphi P}=0$;

\noindent (18) if $(a=4,\ \delta_P=0,\ n\equiv 0\mod 8,\ r\equiv
6\mod 8,\ \delta_{\varphi}=1)$, then $\delta_{\varphi P}=0$.

\medskip

\noindent {\bf Relations near the boundaries $a=0$ and $a=22-r$:}

\noindent (19) if $r=20$, then $n=2^\epsilon p_1^{\alpha_1}\cdots
p_m^{\alpha_m}$, where $\epsilon \le 2$, $p_i$ a prime, $p_i\equiv
1\mod 4$.

\medskip

\begin{figure}
\begin{picture}(200,140)
\put(66,110){\circle{7}} \put(71,108){{\tiny means
$\delta_\varphi=0$}} \put(66,102){\circle*{3}} \put(71,100){{\tiny
means $\delta_\varphi=1$}}

\multiput(8,0)(8,0){20}{\line(0,1){94}}
\multiput(0,8)(0,8){11}{\line(1,0){170}}
\put(0,0){\vector(0,1){100}} \put(0,0){\vector(1,0){180}} \put(
6,-10){{\tiny $1$}} \put( 14,-10){{\tiny $2$}} \put(
22,-10){{\tiny $3$}} \put( 30,-10){{\tiny $4$}} \put(
38,-10){{\tiny $5$}} \put( 46,-10){{\tiny $6$}} \put(
54,-10){{\tiny $7$}} \put( 62,-10){{\tiny $8$}} \put(
70,-10){{\tiny $9$}} \put( 76,-10){{\tiny $10$}} \put(
84,-10){{\tiny $11$}} \put( 92,-10){{\tiny $12$}}
\put(100,-10){{\tiny $13$}} \put(108,-10){{\tiny $14$}}
\put(116,-10){{\tiny $15$}} \put(124,-10){{\tiny $16$}}
\put(132,-10){{\tiny $17$}} \put(140,-10){{\tiny $18$}}
\put(148,-10){{\tiny $19$}} \put(156,-10){{\tiny $20$}}

\put(-8, -1){{\tiny $0$}} \put(-8,  7){{\tiny $1$}} \put(-8,
15){{\tiny $2$}} \put(-8, 23){{\tiny $3$}} \put(-8, 31){{\tiny
$4$}} \put(-8, 39){{\tiny $5$}} \put(-8, 47){{\tiny $6$}} \put(-8,
55){{\tiny $7$}} \put(-8, 63){{\tiny $8$}} \put(-8, 71){{\tiny
$9$}} \put(-10, 79){{\tiny $10$}} \put(-10, 87){{\tiny $11$}}
\put( -2,114){{\footnotesize $a$}} 
\put(186, -2){{\footnotesize $r$}} 

\put( 16, 0){\circle{7}} \put( 80, 0){\circle{7}} \put(144,
0){\circle{7}} \put( 16,16){\circle{7}} \put( 48,16){\circle{7}}
\put( 80,16){\circle{7}} \put(112,16){\circle{7}}
\put(144,16){\circle{7}} \put( 48,32){\circle{7}} \put(
80,32){\circle{7}} \put(112,32){\circle{7}} \put(
80,48){\circle{7}} \put( 112,48){\circle{7}} \put(
80,64){\circle{7}}
\put( 80,80){\circle{7}}    

\put(  8, 8){\circle*{3}} \put( 24, 8){\circle*{3}} \put( 72,
8){\circle*{3}} \put( 88, 8){\circle*{3}} \put(136,
8){\circle*{3}} \put(152, 8){\circle*{3}} \put(
16,16){\circle*{3}} \put( 32,16){\circle*{3}} \put(
64,16){\circle*{3}} \put( 80,16){\circle*{3}} \put(
96,16){\circle*{3}} \put(128,16){\circle*{3}}
\put(144,16){\circle*{3}} \put( 24,24){\circle*{3}} \put(
40,24){\circle*{3}} \put( 56,24){\circle*{3}} \put(
72,24){\circle*{3}} \put( 88,24){\circle*{3}}
\put(104,24){\circle*{3}} \put(120,24){\circle*{3}}
\put(136,24){\circle*{3}} \put( 32,32){\circle*{3}} \put(
48,32){\circle*{3}} \put( 64,32){\circle*{3}} \put(
80,32){\circle*{3}} \put( 96,32){\circle*{3}}
\put(112,32){\circle*{3}} \put(128,32){\circle*{3}} \put(
40,40){\circle*{3}} \put( 56,40){\circle*{3}} \put(
72,40){\circle*{3}} \put( 88,40){\circle*{3}}
\put(104,40){\circle*{3}} \put(120,40){\circle*{3}} \put(
48,48){\circle*{3}} \put( 64,48){\circle*{3}} \put(
80,48){\circle*{3}} \put( 96,48){\circle*{3}}
\put(112,48){\circle*{3}} \put( 56,56){\circle*{3}} \put(
72,56){\circle*{3}} \put( 88,56){\circle*{3}}
\put(104,56){\circle*{3}} \put( 64,64){\circle*{3}} \put(
80,64){\circle*{3}} \put( 96,64){\circle*{3}} \put(
72,72){\circle*{3}} \put( 88,72){\circle*{3}}
\put( 80,80){\circle*{3}}   
\put( 88,88){\circle*{3}} \put( 96,80){\circle*{3}} \put(
104,72){\circle*{3}} \put( 112,64){\circle*{3}} \put(
120,56){\circle*{3}} \put( 128,48){\circle*{3}} \put(
136,40){\circle*{3}} \put( 144,32){\circle*{3}} \put(
144,32){\circle{7}} \put( 152,24){\circle*{3}} \put(
160,16){\circle*{3}}
\end{picture}

\caption{All possible $(r, a, \delta_\varphi)$} \label{Sgraph}
\end{figure}
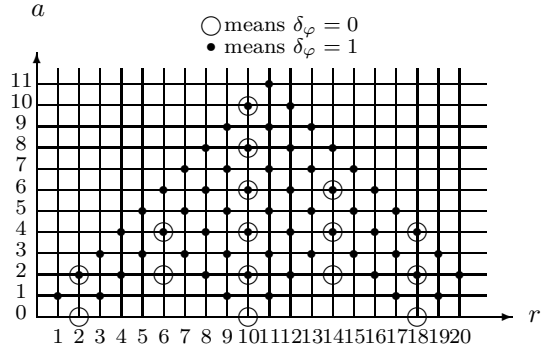

We remark that the invariants $(r,a,\delta_{\varphi})$ satisfying
the conditions 0.(1)---(7) classify, up to isomorphism, all
integral involutions $(L,\varphi)$ satisfying the condition:
$L^\varphi$ is hyperbolic. All these invariants are listed in
Figure \ref{Sgraph}. They classify connected components of moduli
of real K\"ahlerian (and then without a polarization) K3 surfaces.

All polarization conditions I.(1)---(19) depend on $n\mod 8$
except two conditions:  the condition I.(5) depends on $n\mod 16$
and may happen only on the boundary $(r+a=22,\ \delta_\varphi=1)$
(see Figure \ref{Sgraph}); the condition I.(19) depends on prime
decomposition of $n$ and may happen only in one point
$(r,a,\delta_{\varphi })=(20,2,1)$ (see Figure \ref{Sgraph}).
Thus, depending on $n\mod 8$ or $n\mod 16$, all genus invariants
\eqref{genK3} can be easily enumerated similarly to Figure
\ref{Sgraph}. See Figures 33---41 in \cite{NikulinSaito05}.

\medskip

At last, it was shown in \cite{Nikulin79} that ``almost in all cases"
the genus invariants \eqref{genK3} determine the isomorphism class
of an integral polarized K3 involution, and then (by Theorem
\ref{intpolinvtheorem}) the moduli space of real polarized K3
surfaces with this invariants is connected. Exactly we have (see
Theorem 3.3.1 in \cite{Nikulin79}):

{\it The moduli space of real polarized K3 surfaces with the fixed
genus invariants \eqref{genK3} is connected if $r\le 18-\delta_P$.
In particular, it is connected if $r\le 17$.}

\medskip

Looking at Figure \ref{Sgraph}, one can see that there are very
few points $(r,a,\delta_{\varphi})$ which don't satisfy this
condition. But, there are still many.

In the next Section \ref{new} we shall improve and finalize this
connectedness result. When it is not valid, we shall enumerate
connected components of moduli of real polarized K3 surfaces by
additional invariants.

\subsection{The enumeration of connected components of
moduli of real polarized K3 surfaces.} \label{new}
We have the following significant improvement of the old
connectedness result above where the last statement of the theorem
follows from \eqref{realcomponents} and Figure \ref{Sgraph}.

\begin{theorem} The genus invariants \eqref{genK3} determine the 
isomorphism class 
of an integral polarized K3 involution if $r\le 18$.

In particular (by Theorem \ref{intpolinvtheorem}) the moduli space
of real polarized K3 surfaces with the fixed genus invariants
\eqref{genK3} is connected if $r\le 18$.

In particular, the moduli space of real polarized K3 surfaces
$(X,P^\prime)$ with the fixed genus invariants \eqref{genK3} is
connected if $X(\br)$ is different from $T_1\amalg (T_0)^9$ (i. e.
$(r,a)=(19,1)$), and $(T_0)^9$ (i. e. $(r,a)=(19,3)$), and
$(T_0)^{10}$ (i. e. $(r,a)=(20,2)$). \label{theoremuniqueness}
\end{theorem}

\Proof Fix the genus invariants \eqref{genK3} satisfying $r\le
18$. The invariants determine the genus of the lattices
$L^\varphi$ and $L_{\varphi,P}=(L^\varphi\oplus \bz
P)^\perp_L=(P)^\perp_{L_{\varphi}}$. We denote by $q_{S}$ the
discriminant quadratic form of an even lattice $S$ (see
\cite{Nikulin79}). The statement follows if the lattices
$L^\varphi$ and $L_{\varphi,P}$ are unique in their genus, and the
canonical homomorphisms
$$
O(L^{\varphi})\to O(q_{L^{\varphi}}),\ \ \ O(L_{\varphi,P})\to
O(q_{L_{\varphi,P}})
$$
are epimorphic. See the proof of Theorem 3.3.1 in \cite{Nikulin79}
(or look at Remark 1.6.2 of \cite{Nikulin83}). Here $O$ denotes
the full orthogonal group.

The lattice $L^\varphi$ is a hyperbolic 2-elementary (i.
e. $(L^\phi)^\ast/L^\phi\cong (\bz/2\bz)^a$ is a 2-elementary
Abelian group) lattice. For this lattices the statement is always
valid (Theorems 3.6.2, 3.6.3 in \cite{Nikulin79}).

The lattice $M=L_{\varphi,P}$ is an indefinite lattice of the rank
at least $3$ since $\rk M=22-r-1$ and $r\le 18$. Its discriminant
group $A_M=M^\ast/M$ and the discriminant form on that group are
calculated in \cite{Nikulin79} (see also \cite{Nikulin83}). Its
$p$-component is a cyclic group $(\bz/p^{k_p})$ for any odd prime
$p$. Its $2$-component is isomorphic to $(\bz/2\bz)^{k_1}\oplus
(\bz/2^k\bz)$ where $k\ge 0$.

Assume that $\rk M\le 17$. Then $\rk M=21-r\ge 4$, and the
required statement for $M=L_{\varphi,P}$ follows from Theorem
1.14.2 in \cite{Nikulin79} (or Theorem 1.2$^\prime$ in
\cite{Nikulin79.1}).

Assume that $r=18$. Then $\rk M=3$. If the 2-component of $A_M$ is
a cyclic group, then the statement again follows from Theorem
1.14.2 in \cite{Nikulin79} (or Theorem 1.2$^\prime$ in
\cite{Nikulin79.1}). For instance, it is valid if the invariant
$\delta_P=0$.

Assume that $r=18$, but the 2-component of $A_M$ is not cyclic.
Since $M$ is even, it follows that the 2-component of $A_M$ is
isomorphic to $(\bz/2\bz)^2\oplus (\bz/2^k\bz)$ where $k\ge 1$. In
this case, the statement follows from the general results
announced by Miranda and Morrison in \cite{MirandaMorrison85},
\cite{MirandaMorrison86}. The proofs are contained in James
\cite{James96}. Remark that we need a very particular case of the
general results of Miranda---Morrison and James which is a little
bit stronger than Theorem 1.14.2 in \cite{Nikulin79}.

It follows the theorem. \QED

\medskip

Now let us assume that $r\ge 19$. Then the lattice
$M=L_{\varphi,P}$ has the rank 2 or 1, and the statement of
Theorem \ref{theoremuniqueness} surely is not valid (see below). We must
introduce some additional invariants.

By Theorem \ref{geninvtheorem} (see Figure \ref{Sgraph}), if $r\ge
19$, then $(r,a,\delta_\varphi)=(19,1,1)$, $(19,3,1)$ or
$(20,2,1)$. Assuming one of this cases, we apply the following
general construction which uses the technique of discriminant
forms from \cite{Nikulin79}.

We fix the invariants \eqref{genK3} with the $(r,a,\delta_\varphi)$ where
$r\ge 19$. The invariants $(r,a,\delta_\varphi)$ define the isomorphism
class of the lattice $L^\varphi$. The additional invariants
$(n,\delta_P,\delta_{\varphi P})$ define the lattice
\begin{equation}
L^{\varphi,P}=[L^\varphi,P]_{pr} \label{LplusP}
\end{equation}
which is a primitive sublattice in $L$ generated by $L^\varphi$
and $P$. We denote by
\begin{equation} q=q_{L^{\varphi,P}}
\label{q}
\end{equation}
the discriminant quadratic form of the lattice $L^{\varphi,P}$. We
introduce the group 
$$
O(L^\varphi)_P=\{f\in O(L^{\varphi,P})\ |\ f(P)=P\}\subset
O(L^\varphi),
$$
and calculate the image of the natural homomorphism
\begin{equation}
O(L^\varphi)_P\to \overline{O(L^\varphi)_P}\subset O(q),
\label{Oplusq}
\end{equation}
which we denote as $\overline{O(L^\varphi)_P}$. The lattice
$L_{\varphi,P}$ is orthogonal to $L^{\varphi,P}$ in the even
unimodular lattice $L$. Then it has the signature $(1,20-r)$ and
the discriminant quadratic form which is isomorphic to $(-q)$.
This fixes the genus of the lattice $L^{\varphi,P}$. Any lattice
$L^{\varphi,P}$ with these invariants can be taken, and its class
of isomorphism is the main additional invariant of the integral
polarized K3 involution.

Since $L^{\varphi,P}$ and $L_{\varphi,P}$ are orthogonal
complements to each other in the even unimodular lattice $L$, it
defines the canonical isomorphism
\begin{equation}
\tau:q_{L_{\varphi,P}}\cong -q
\label{tau}
\end{equation}
of the finite quadratic forms. Any of them can be taken, and it
defines another invariant of an integral polarized K3 involutions.
Any such a pair
\begin{equation}
(L_{\varphi,P}\,,\ \tau)
\label{addinvariants}
\end{equation}
can be taken, and it defines the isomorphism class of an integral
polarized K3 involution with the given genus invariants. More
exactly,
\begin{equation}
L=[L^{\varphi,P}\oplus L_{\varphi,P},\ \{x^\ast_+\oplus x^\ast_-\
|\ x^\ast_+\in (L^{\varphi,P})^\ast,\ x^\ast_-\in
(L_{\varphi,P})^\ast\ \text{and}\ \tau (x^\ast_- +
L_{\varphi,P})=x^\ast_+ + L^{\varphi,P} \}]
\label{constLphi}
\end{equation}
with the action of the involution $\varphi$ which is +1 on
$L^\varphi$, and $-1$ on $P$ and $L_{\varphi,P}$. Other speaking
$$
L/\left(L^{\varphi,P}\oplus L_{\varphi,P}\right)\subset
(L^{\varphi,P})^\ast/L^{\varphi,P}\oplus (L_{\varphi,P})^\ast/L_{\varphi,P}
$$
is the graph of the isomorphism $\tau$.

Two such pairs $(L_{\varphi,P}\,,\ \tau)$ and
$(\widetilde{L_{\varphi,P}}\,,\ \widetilde{\tau })$ define
isomorphic integral polarized K3 involutions if and only if there
exists an isomorphism $f:L_{\varphi,P}\to
\widetilde{L_{\varphi,P}}$ of lattices such that for the induced
isomorphism $\overline{f}:q_{L_{\varphi,P}}\to
q_{\widetilde{L_{\varphi,P}}}$ of their discriminant forms, one
has
\begin{equation}
\widetilde{\tau}\overline{f}=g\tau
\end{equation}
where $g\in \overline{O(L^{\varphi })_P}$.

Thus, for the fixed isomorphism class $L_{\varphi,P}$, the number
of classes of pairs is equal to the number of double cosets
\begin{equation}
\overline{O(L_{\varphi,P})}\backslash
O(q)/\overline{O(L^\varphi)_P} \label{doublecosets}
\end{equation}
where $\overline{O(L_{\varphi,P})}$ is the image of
$O(L_{\varphi,P})$ in $O(q)$ by using $\tau$.

Thus, the actual number of classes of integral polarized K3
involutions and connected components of moduli of real polarized
K3 surfaces is equal to
\begin{equation} \sum_{classes
L_{\varphi,P}\ of\ genus\
(1,20-r,-q)}{\sharp(\overline{O(L_{\varphi,P})}\backslash
O(q)/\overline{O(L^\varphi)_P})}. \label{numbofclasses}
\end{equation}

We remark that
\begin{equation}
O(q)=\prod_pO(q_p),
\label{decOq}
\end{equation}
where $q_p$ is the non-trivial $p$-component of $q$ for a prime
$p$. The $p$-component $q_p$ is a quadratic form on a cyclic group
if $p|n$ is odd, and then $O(q_p)=\{\pm 1_p\}$ where $-1_p$ means
minus identity on the $p$-component $q_p$. The same is valid for
$q_2$ if the $2$-component is also cyclic. We also remark that
$O(L^\varphi)_P$ acts only on $2$-component of $q$, and it is
trivial almost in all cases. Thus, we have
\begin{equation}
O(q)=O(q_2)\prod_{\text{odd\ }p|n}\{\pm 1_p\},\ \
\overline{O(L^\varphi)}_P\subset O(q_2). \label{decoq2}
\end{equation}
Below we calculate these
invariants for each of three cases $(r,a,\delta_\varphi)$ with
$r\ge 19$.

\medskip

We use notations: $\langle A\rangle$ denote a lattice defined by a
symmetric integral matrix $A$. We denote $U=
\left\langle\begin{array}{cc}
0 & 1\\
1 & 0
\end{array}\right\rangle$. It is a unique class of an even unimodular 
lattice of the signature $(1,1)$. 
The lattice $E_8$ is a negative definite even unimodular lattice of the
rank 8. It is defined by the corresponding diagram $E_8$. For a
lattice $K$, we denote by $K(m)$ a lattice obtained from the
lattice $K$ by multiplication of the form of $K$ by $m\in \bq$. By
$\oplus$ we denote the orthogonal sum of lattices.

\subsubsection{The case $(r,a,\delta_\varphi)=(19,1,1)$.}
\label{subsubsec1911} In this case $L^\varphi=U\oplus 2E_8\oplus
\langle -2 \rangle$.

{\bf The case $\delta_P=1$.} Then $L^{\varphi,P}=L^\varphi\oplus
\bz P=U\oplus 2E_8\oplus \langle -2 \rangle\oplus \langle
n\rangle$. Since $U$ and $E_8$ are unimodular, the discriminant
form
\begin{equation}
q=q_{L^{\varphi,P}}=q_{\langle -2\rangle \oplus \langle n\rangle}.
\end{equation}
The group $\overline{O(L^{\varphi,P})_P}$ is trivial since it acts
only on the orthogonal part $q_{\langle -2 \rangle}$ which has the
group of order 2. The hyperbolic lattice $\langle 2 \rangle \oplus
\langle -n \rangle$ has the discriminant quadratic form which is isomorphic to
$-q$. Thus, $L^{\varphi, P}$ is any lattice having the same genus. Here
$n>0$ can be any even number, the invariant $\delta_{\varphi P}=1$
(see Theorem \ref{geninvtheorem}). We naturally identify
$-q=-q_{\langle -2\rangle \oplus \langle n\rangle}=
q_{\langle 2\rangle \oplus \langle -n\rangle}$.  Then an isomorphism
$\tau:q_{L_{\varphi,P}}\to -q$ is identified with the isomorphism
$\tau:q_{L_{\varphi,P}}\to
q_{\langle 2 \rangle} \oplus q_{\langle -n \rangle}$.

Thus, we finally obtain the result.

\begin{theorem}
\label{theorem1911} If
$(r,a,\delta_{\varphi},\delta_P)=(19,1,1,1)$ (equivalently,
$X(\br)=T_1\amalg (T_0)^8$ and $\delta_P=1$), then
$\delta_{\varphi P}=1$, $n$ any positive even integer, $k$ any
positive integer, and the isomorphism classes of integral
polarized K3 involutions (equivalently connected components of
moduli of real polarized K3 surfaces) with these invariants are in
one to one correspondence with pairs
$$
(L_{\varphi,P}, \tau)
$$
where $L_{\varphi,P}$ is any lattice of the genus $\langle 2
\rangle\oplus \langle -n \rangle$ and $\tau:q_{L_{\varphi,P}}\to
q_{\langle 2 \rangle} \oplus q_{\langle -n \rangle}$ is any
isomorphism of quadratic forms considered up to the natural action
of $O(L_{\varphi,P})$.

Thus, the number of connected components of moduli is equal to
$$
\sum_{L_{\varphi,P}\ of\ genus\ \langle 2 \rangle \oplus \langle
-n\rangle}{\sharp(O(q_{L_{\varphi,P}})/\overline{O(L_{\varphi,P})}
)}.
$$
\end{theorem}

\medskip

For our further study, the following definition is very important.
{\it A connected component of moduli} from Theorem \ref{theorem1911}
{\it is called standard} if $L_{\varphi,P}=\langle 2 \rangle \oplus \langle -n
\rangle$ (it defines the identity of the corresponding Gauss class group
of binary lattices, e. g. see \cite{Cassels78}) and $\tau$ is induced by
the identity isomorphism of the lattice
$\langle 2\rangle \oplus \langle -n \rangle$.
Any other connected component of moduli with the lattice $L_{\varphi,P}\cong
\langle 2 \rangle \oplus \langle -n \rangle$ will be different
from the standard one by an automorphism of
$O(q_{\langle 2 \rangle}\oplus q_{\langle -n\rangle})$ (up
to $O(\langle 2 \rangle \oplus \langle -n \rangle)$),
and it can be labelled by the automorphism. If the automorphism belongs 
to the 2-component 
$O(q_{\langle 2 \rangle}\oplus q_{{\langle -n\rangle}_2})$ of
the automorphism group, we say that this component is
{\it different from the standard only over $2$}. Of course, 
we use the same names  
for real polarized K3 surfaces from these connected components of moduli.
Similar definitions we use in all cases below when we introduce
the standard connected component of moduli.

We remark that the 2-component $O(q_{\langle 2 \rangle}\oplus
q_{{\langle -n \rangle}_2})$ is trivial if $n\equiv 2\mod 8$, it
is $(\bz/2)^2$ if $n\equiv 0\mod 16$, and it is $\bz/2$ in the
remaining cases.

Thus, the number of connected components of moduli which
are different from the standard only over $2$ is one if $n
\equiv 2\mod 8$, it is at most four if $n\equiv 0\mod 16$, and it is
at most two in the remaining cases.

For example, Theorem \ref{theorem1911} gives only one connected component of
moduli which is the standard one if
$n=2$, $4$, $6$ or $8$.

\medskip

{\bf The case $\delta_P=0$.} Since
$(r,a,\delta_{\varphi})=(19,1,1)$ and $\delta_P=0$, by Theorem
\ref{geninvtheorem} (see relations I.(2),(9),(19)) we have
$\delta_{\varphi P}=0$, $P^2=n\equiv 2\mod 8$, and the lattice
$\bz P =\langle n \rangle$.

We again write $L^\varphi=U\oplus 2E_8\oplus \langle -2 \rangle$.
We denote by $e$ the generator of the summand $\langle -2 \rangle$
with $e^2=-2$. Since $\delta_P=0$ and $a=1$,
$$
L^{\varphi,P}=U\oplus 2E_8\oplus \left(\langle -2 \rangle \oplus
\langle n \rangle \right)(1/2,1/2)
$$
where $\langle -2 \rangle \oplus \langle n \rangle(1/2,1/2)\supset
\langle -2 \rangle \oplus \langle n \rangle$ is the overlattice of
the index 2 defined by
\begin{equation}
\langle -2 \rangle \oplus \langle n \rangle (1/2,1/2)=[e,
e/2+P/2]\,= \left\langle\begin{array}{cc}
-2 & -1\\
-1 & \frac{n-2}{4}\\
\end{array}\right\rangle.
 \label{latt2}
\end{equation}
Since $U$ and $E_8$ are unimodular and the lattice \eqref{latt2}
is unimodular over 2, it follows
\begin{equation}
q_{L^{\varphi,P}}=q\left(\left\langle\begin{array}{cc}
-2 & -1\\
-1 & \frac{n-2}{4}\\
\end{array}\right\rangle\right)=\oplus_{\text{odd}\ p\vert n}q_{{\langle
n\rangle}_p}.
\end{equation}
Here we denote by $q(S)$ (or $b(S))$ the discriminant quadratic (or bilinear)
form of a lattice $S$.

Then the group $\overline{O(L^\varphi)_P}$ is again trivial, and
we obtain similarly to the previous case

\begin{theorem}
\label{theorem1910} If
$(r,a,\delta_{\varphi},\delta_P)=(19,1,1,0)$ (equivalently,
$X(\br)=T_1\amalg (T_0)^8$ and $\delta_P=0$), then
$\delta_{\varphi P}=0$, $n$ any positive integer such that
$n\equiv 2\mod 8$, $k$ any positive integer, and the isomorphism
classes of integral polarized K3 involutions (equivalently
connected components of moduli of real polarized K3 surfaces) with
these invariants are in one to one correspondence with pairs
$$
(L_{\varphi,P}, \tau)
$$
where $L_{\varphi,P}$ is any lattice of the genus $
\left\langle\begin{array}{cc}
2 & 1\\
1 & \frac{2-n}{4}\\
\end{array}\right\rangle
$, and
$$
\tau:q_{L_{\varphi,P}}\to
q\left(\left\langle\begin{array}{cc}
2 & 1\\
1 & \frac{2-n}{4}\\
\end{array}\right\rangle\right)=
\oplus_{odd\ p|n}{q_{{\langle
-n \rangle}_p}}
$$
is any isomorphism of quadratic forms considered
up to the natural action of $O(L_{\varphi,P})$.

Thus, the number of connected components of moduli is equal to
$$
\sum_{L_{\varphi,P}^{}\ of\ genus\ \left\langle\begin{array}{cc}
2 & 1\\
1 & \frac{2-n}{4}\\
\end{array}\right\rangle}
{\sharp(O(q_{L_{\varphi,P}})/\overline{O(L_{\varphi,P})} )}.
$$
\end{theorem}

\medskip

We define the {\it standard connected component of moduli} when
$L_{\varphi,P}=$ $\left\langle\begin{array}{cc}
2 & 1\\
1 & \frac{2-n}{4}\\
\end{array}\right\rangle
$ and $\tau$ is induced by the identity isomorphism of the lattice.

In this case, any connected component of moduli which is different
from the standard only over $2$ coincides with the standard one.

For example, Theorem \ref{theorem1910} gives only one connected
component of moduli which is the standard one if $n=2$.

\subsubsection{The case $(r,a,\delta_\varphi)=(19,3,1)$.}
\label{subsubsec193} By Theorem \ref{geninvtheorem}, relation
I,(5), we have: if $\delta_{\varphi P}=0$, then $n\equiv 2\mod
16$. It is the only relation between genus invariants
\eqref{genK3} which we have in this case.

We apply a modification of the general construction above.

Since $(r,a,\delta_\varphi)=(19,3,1)$, the lattice $L_\varphi$
where $\varphi=-1$ is an even 2-elementary lattice of signature
$(2,1)$ and with $a=3$. It follows that $L_\varphi(1/2)$ is an
even unimodular lattice of signature $(2,1)$. This lattice is odd,
and we have
$$
L_\varphi(1/2)\cong \langle 1\rangle \oplus \langle 1 \rangle
\oplus \langle -1 \rangle .
$$
by classification of unimodular indefinite lattices.

The element $P \in L_\varphi(1/2)$ is then a primitive element
with $P^2=n/2$. Since both lattices $L^\varphi$ and $L_\varphi$
are unique in their genus and for both of them the canonical
homomorphisms
$$
O(L^\varphi)\to O(q_{L^\varphi}),\ \ O(L_\varphi)\to
O(q_{L_\varphi})
$$
are epimorphic, any automorphism of $L_\varphi$ can be extended to
an automorphism of $L$ which commutes with the involution
$\varphi$.

Thus, two integral polarized involutions are isomorphic if and
only if the corresponding to them elements $P\in L_\varphi(1/2)$
are conjugate by $O(L_{\varphi})$. The invariant $\delta_{\varphi P}$
has then the following meaning. We have
$P^\perp_{L_{\varphi}(1/2)}=L_{\varphi,P}(1/2)$ is odd if
$\delta_{\varphi P}=1$, and it is even if $\delta_{\varphi P}=0$.
The last case is possible only if $n\equiv 2\mod 16$.

The sublattice $\bz P\subset L_\varphi (1/2)$ has the discriminant
bilinear form $b_{\langle n/2\rangle}$. Then the lattice
$L_{\varphi,P}(1/2)$ has the discriminant bilinear form
$-b_{\langle n/2 \rangle}$ and signature $(1,1)$.

If $\delta_{\varphi P}=1$, this lattice is odd, and has the genus
of the lattice
$$
\langle 1 \rangle\oplus \langle -n/2 \rangle .
$$
Thus $L_{\varphi,P}$ has the genus of the lattice
$$
\langle 2 \rangle \oplus \langle -n \rangle .
$$

If $\delta_{\varphi P}=0$, we should change the lattice $\langle 1
\rangle\oplus \langle -n/2 \rangle$ to make it even, but having
the same bilinear discriminant form. Let us denote by $e_1$ and
$e_2$ the generators of $\langle 1\rangle$ and $\langle -n/2
\rangle$ respectively. The new lattice will be
$[2e_1,2e_2,(e_1+e_2)/2]$. This changes the lattice only over $2$.
The new lattice has the basis $\{2e_1,(e_1+e_2)/2\}$, and has the
matrix
$$
\left\langle\begin{array}{cc}
4 & 1\\
1 & \frac{2-n}{8}\\
\end{array}\right\rangle.
$$
It is even and unimodular over 2 if $n\equiv 2\mod 16$. Thus, the
lattice $L_{\varphi,P}$ has the genus of the lattice
$$
\left\langle\begin{array}{cc}
8 & 2\\
2 & \frac{2-n}{4}\\
\end{array}\right\rangle.
$$
Finally, we obtain

\begin{theorem}
\label{theorem1931} If $(r,a,\delta_{\varphi},\delta_{\varphi
P})=(19,3,1,1)$ (equivalently,
$X(\br)=(T_0)^9$ and $X(\br)\not\sim P\mod 2$ in $H_2(X(\bc),\bz)$),
then $\delta_P=0$, $n$ any positive even integer,
$k$ any positive integer, and the isomorphism classes of integral
polarized K3 involutions (equivalently connected components of
moduli of real polarized K3 surfaces) with these invariants are in
one to one correspondence with pairs
$$
(L_{\varphi,P}, \tau)
$$
where $L_{\varphi,P}$ is any lattice of the genus $\langle 2
\rangle\oplus \langle -n \rangle$, and
$$
\tau:b_{L_{\varphi,P}(1/2)}\to b_{\langle 1\rangle \oplus  
\langle -n/2 \rangle}=
b_{\langle -n/2 \rangle}
$$
is any isomorphism of bilinear forms considered up to the natural
action of $O(L_{\varphi,P})$.

Thus, the number of connected components of moduli is equal to
$$
\sum_{L_{\varphi,P}\ of\ genus\ \langle 2 \rangle \oplus \langle
-n\rangle
}{\sharp(O(b_{L_{\varphi,P}(1/2)})/\overline{O(L_{\varphi,P})} )}.
$$
\end{theorem}

\medskip

We define the {\it standard connected component of moduli} when
$L_{\varphi,P}=\langle 2 \rangle\oplus \langle -n \rangle$ and $\tau$ is
induced by the identity isomorphism of the lattice.

We remark that the 2-component $O(b_{{\langle -n/2 \rangle}_2})$
is trivial if $n\equiv 2\mod 4$ or $n\equiv 4\mod 8$, it is
$\bz/2$ if $n\equiv 8\mod 16$, and it is $\bz/2\times \bz/2$ if
$n\equiv 0\mod 16$.

Thus, the number of connected components of moduli which are different from
the standard only over $2$ is one if $n\equiv 2\mod 4$ or $n\equiv 4\mod 8$;
it is at most two if $n\equiv 8\mod 16$; it is at most four if
$n\equiv 0\mod 16$.

For example, Theorem \ref{theorem1931} gives only one connected component
of moduli which is the standard one if $n=2$, $4$, $6$ or $8$.

\medskip

\begin{theorem}
\label{theorem1930} If $(r,a,\delta_{\varphi},\delta_{\varphi
P})=(19,3,1,0)$  (equivalently,
$X(\br)=(T_0)^9$ and $X(\br)\sim P\mod 2$ in $H_2(X(\bc),\bz)$),
then $\delta_P=0$, $n$ any positive integer such
that $n\equiv 2\mod 16$, $k$ any positive integer, and the
isomorphism classes of integral polarized K3 involutions
(equivalently connected components of moduli of real polarized K3
surfaces) with these invariants are in one to one correspondence
with pairs
$$
(L_{\varphi,P}, \tau)
$$
where $L_{\varphi,P}$ is any lattice of the genus
$$
\left\langle\begin{array}{cc}
8 & 2\\
2 & \frac{2-n}{4}\\
\end{array}\right\rangle,
$$ and
$$
\tau:b_{L_{\varphi,P}(1/2)}\to b\left(\left\langle\begin{array}{cc}
4 & 1\\
1 & \frac{2-n}{8}\\
\end{array}\right\rangle\right)=
b_{\langle -n/2 \rangle}
$$
is any isomorphism of bilinear forms considered up to the natural
action of $O(L_{\varphi,P})$.

Thus, the number of connected components of moduli is equal to
$$
\sum_{L_{\varphi,P}\ of\ genus\  \left\langle\begin{array}{cc}
8 & 2\\
2 & \frac{2-n}{4}\\
\end{array}\right\rangle
}{\sharp(O(b_{L_{\varphi,P}(1/2)})/\overline{O(L_{\varphi,P})} )}.
$$
\end{theorem}

\medskip

We define the {\it standard connected component of moduli} when
$L_{\varphi,P}=\left\langle\begin{array}{cc}
8 & 2\\
2 & \frac{2-n}{4}\\
\end{array}\right\rangle$ and $\tau$ is induced by the identity isomorphism
of the lattice.

Since $n\equiv 2 \mod 4$ in this case, any connected component of moduli
which is different from the standard only over 2 is the standard one.

For example, Theorem \ref{theorem1930} gives only one connected component
of moduli which is the standard one if $n=2$.

\medskip

\subsubsection{The case $(r,a,\delta_\varphi)=(20,2,1)$.}
\label{subsubsec202} By Theorem \ref{geninvtheorem}, relations
I.(4),(11),(15),(19), we have: $\delta_P=0$; $n=2^\epsilon
p_1^{\alpha_1}\cdots p_m^{\alpha_m}$, where $1\le \epsilon \le 2$,
$p_i$ a prime, $p_i\equiv 1\mod 4$; $\delta_{\varphi P}=1$ if
$n\equiv 2\mod 4$ (i. e. $\epsilon=1$), and $\delta_{\varphi P}=0$
if $n\equiv 0\mod 4$ (i. e. $\epsilon =2$).

This case is very similar to the previous one. Since
$(r,a,\delta_\varphi)=(20,2,1)$, the lattice $L_\varphi$ is a
positive definite even 2-elementary lattice of the rank $2$ and
with $a=2$. It follows that $L_\varphi(1/2)$ is an even unimodular
positive definite lattice of the rank $2$. This lattice is odd and
$L_\varphi(1/2)\cong \langle 1\rangle \oplus \langle 1 \rangle$ by
classification of unimodular lattices of a small rank.

We have the same results as for previous case.

Two integral polarized involutions are isomorphic if and only if
the corresponding to them elements $P\in L_\varphi(1/2)$ are
conjugate by $O(L_{\varphi})$. We have
$P^\perp_{L_{\varphi}(1/2)}=L_{\varphi,P}(1/2)$ is odd if
$\delta_{\varphi P}=1$, and it is even if $\delta_{\varphi P}=0$.

The sublattice $\bz P\subset L_\varphi (1/2)$ has the discriminant
bilinear form $b_{\langle n/2\rangle}$. Then the lattice
$L_{\varphi,P}(1/2)$ has the discriminant bilinear form
$-b_{\langle n/2 \rangle}$, and it is positive definite of the
rank one. It follows that $L_{\varphi,P}\cong \langle n\rangle$,
and the bilinear form $b_{\langle n/2 \rangle}$ is isomorphic to
$-b_{\langle n/2\rangle}$. This is equivalent to the condition on
$n$ above.

Thus, we obtain

\begin{theorem}
\label{theorem202} If $(r,a,\delta_{\varphi})=(20,2,1)$ (equivalently,
$X(\br)=(T_0)^{10}$), then
$\delta_P=0$; $\delta_{\varphi P}=1$ (equivalently, $X(\br)\not\sim P\mod 2$
in $H_2(X(\bc),\bz)$)
if $n\equiv 2\mod 4$, and
$\delta_{\varphi P}=0$ (equivalently, $X(\br)\sim P\mod 2$
in $H_2(X(\bc),\bz)$) if $n\equiv 0\mod 4$; $n=2^\epsilon
p_1^{\alpha_1}\cdots p_m^{\alpha_m}$, where $\epsilon=1$ or $2$,
$p_i$ a prime, $p_i\equiv 1\mod 4$; $k$ any positive integer.
Isomorphism classes of integral polarized K3 involutions
(equivalently connected components of moduli of real polarized K3
surfaces) with these invariants are in one to one correspondence
with isomorphisms $\tau:b_{\langle n/2 \rangle}\to -b_{\langle
n/2\rangle}$ of finite bilinear forms considered up to $\pm 1$.
(The lattice $L_{\varphi,P}=\langle n\rangle$. ) Thus, the number
of connected components of moduli is equal to $2^{\max\{0,m-1\}}$
where $m$ is the number of different odd prime divisors of $n$.
\end{theorem}

In this case, we don't have a notion of a standard connected
component of moduli.

\section{Deformations of real hyper-elliptically polarized
K3 surfaces.} \label{sechyppolK3}

Here we apply results above to classify real polarized K3 surfaces
which are deformations of real hyper-elliptically polarized K3
surfaces. This question had been studying in
\cite{NikulinSaito05}. Using results of Section \ref{secmoduliK3},
we will be able to finalize these results.

First let us formulate the problem exactly. Let $(X,\wP)$ be a K3
surface with an ample divisor class $\wP$ (not necessarily very
ample) such that the linear system $|\wP|$ does not have fixed
components. Let $P$ be the corresponding primitive element of the
Picard lattice of $X$ such that $\wP=mP$ where $m\in \bn$. As
above, $n=P^2$ is the primitive degree. A pair $(X,\wP)$ is called {\it a
hyper-elliptically polarized K3 surface} if the complete linear
system $|\wP|$ does not give an embedding of $X$ into a projective
space.

Let $(X,\wP)$ be a hyper-elliptically polarized K3 surface. By the
results of Saint-Donat \cite{Saint-Donat74}, a general curve of
the linear system $|\wP|$ is hyper-elliptic in this case, and the
linear system $|\wP|$ gives a double covering $|\wP|:X\to Y\subset
\bp^N$ onto a rational surface $Y$ where
$N=\dim|\wP|=\wP^2/2+1=m^2n/2+1$. Taking $P^\prime=kP$ where $k>m$
a sufficiently large (by \cite{Saint-Donat74}, it is enough to
take $k\ge 3$), we obtain a polarized K3 surface $(X,kP)$ (i. e.
$kP$ is very ample). It is naturally also called {\it
hyper-elliptically polarized}. Then any polarized K3 surface in
the same connected component of moduli as $(X,kP)$ can be
considered as a deformation of the hyper-elliptically polarized K3
surface $(X,kP)$ (or just of the $(X,\wP)$).

It is easy to see that any complex polarized K3 surface is such a
deformation because the moduli space of complex polarized K3
surfaces $(X,kP)$ of the degree $n=P^2>0$ is connected by Global
Torelli Theorem for K3 surfaces \cite{PS71}. It is empty if $n=2$
and $k=1$ or $2$).

We ask similar question for real polarized K3 surfaces: {\it When
a real polarized K3 surface is a deformation a hyper-elliptically
polarized K3 surface?} It could be important because reduces some
questions to hyper-elliptically polarized K3 surfaces which are
much simpler and can be constructed explicitly.

We can reformulate this question as follows. Let us denote by
$\M_{n,k}$ the moduli space of real polarized K3 surfaces $(X,kP)$
where $P$ is primitive, $P^2=n$ and $k\in \bn$ (i. e. $kP$ is very
ample). When $n=2$, we assume that $k\ge 3$, since $\M_{2,1}$ and
$\M_{2,2}$ are empty. We have the obvious embedding of moduli
spaces
\begin{equation}
\M_{n,k_1}\subset \M_{n,k_2}\ if\ k_1\le k_2
\label{concompk1k2}
\end{equation}
and the induced map for the sets of their connected components of
moduli. By Theorem \ref{intpolinvtheorem}, \eqref{concompk1k2}
gives an isomorphism on the sets of connected components. The
connected components of $\M_{n,k}$ are the same for any $k$.
Moreover, it is known that difference in $\M_{n,k}$ for different
$k$ is only in codimension $\ge 1$ (the corresponding K3 surfaces
have Picard number $\ge 2$). Thus, even when $\M_{n,k}$ does not
have hyper-elliptically polarized K3 surfaces, we still can
consider these K3 surfaces as deformations of hyper-elliptically
polarized K3 surfaces if the same connected component of moduli of
$\M_{n,k_2}$ for $k_2>k$ has hyper-elliptically polarized K3
surfaces. Thus, our question does not depend on $k$ and can be
formulated as follows:

{\it Which connected components of $\M_{n,k}$ for $k\gg 0$ contain
hyper-elliptically polarized K3 surfaces $(X,kP)$, i. e. $P$ is
primitive and the linear system $|\wP|$ is hyper-elliptic for
$\wP=mP$ where $m=1$ or $2$?}

This question had been studying in (\cite{NikulinSaito05}, Sect. 8).
Using results of Section \ref{secmoduliK3}, here we finalize these
results.

The exposition in (\cite{NikulinSaito05}, Sect. 8) was very short, and first
we clarify the general considerations in \cite{NikulinSaito05}.

By Theorem 5.2 in \cite{Saint-Donat74}, for an ample $\wP=mP$, the
linear system $|\wP|$ is hyper-elliptic in the following and only
the following cases:

\medskip

\noindent (i) $n=2$, $m=1$ or $2$. Then $Y=\bp^2$, $|\wP|:X\to
Y=\bp^2$ is a double covering ramified in a degree 6 curve, 
$P$ is the preimage of a line in $\bp^2$.

\noindent \noindent (ii) $n\ge 4$, there exists an elliptic curve
$C$ on $X$ such that $\wP\cdot C=2$.

\medskip

It follows that for $n=2$ all polarized K3 surfaces are
hyper-elliptic, and then they are deformations of hyper-elliptic
ones.

Let us assume that $n\ge 4$. Since $C$ is an elliptic curve,
$C^2=0$. Thus $\wP$ and $C$ have the Gram matrix $
2\left(\begin{array}{cc}
n/2 & 1\\
1 & 0\\
\end{array}\right)
$ where the matrix $ \left(\begin{array}{cc}
n/2 & 1\\
1 & 0\\
\end{array}\right)
$ has the determinant $-1$, and it is unimodular. It follows that
the 2-dimensional primitive sublattice $S\subset S_X$ generated by
$\wP$ and $C$ in the Picard lattice $S_X$ of $X$ is a 2-elementary
lattice, i. e. $S^\ast/S\cong (\bz/2\bz)^a$ where $a\le 2$.

Since $P$ is ample, $P^\perp_{S_X}$ has no elements with square
$-2$. It follows that $S^\perp_{S_X}$ also has no elements with
square $-2$. By Global Torelli Theorem for K3 surfaces
\cite{PS71}, there exists an involution $\tau$ of $X$ which is
identity on $S$ and which is $-1$ on the orthogonal complement
$S^\perp$ in $H_2(X(\bc),\bz)$. This involution is non-symplectic,
$Y=X/\{1,\tau\}$ is a rational surface and the quotient map
$\pi:X\to Y$ is a double covering of $Y$ ramified in a
non-singular curve $A\in |-2K_Y|$). See \cite{AlexeevNikulin88}
for details. Depending on the isomorphism
class of $S$, we obtain 5 cases (see the general classification in
\cite{AlexeevNikulin88}, or see \cite{NikulinSaito05}):

{\bf Case $\bff_1$:} $S\cong \langle 2 \rangle\oplus \langle -2
\rangle$. Then $Y=X/\{1,\tau\}=\bff_1$ is a blow-up of $\bp^2$ in
one point. Denoting by $c$ the class of $C$, we obtain that the
elliptic pencil $|c|$ is the preimage of the rational pencil on
$\bff_1$. We denote by $h$ the class of the preimage of a line $l$
in $\bp^2$ and by $e$ the class of the preimage of the exceptional
section $s$ of $\bff_1$. Thus, we have $h^2=2$, $e^2=-2$ and
$h\cdot e=0$. Then $c=h-e$. Since $\wP\cdot c=2$ and $\wP\cdot
e>0$, it follows that $\wP=P=n_1c+e=n_1h+(1-n_1)e$ where
$n=P^2=4n_1-2\equiv 2 \mod 4$ and $n\ge 6$. By Riemann-Roch
Theorem on $\bff_1$, it is easy to see that the linear system
$|P|$ is the preimage of the linear system from $\bff_1$. It
follows that the map $|\wP|:X\to Y\subset \bp^N$ is the quotient map.

{\bf Case $\bff_4$:} $S\cong U$. Then $Y=X/\{1,\tau\}=\bff_4$ where
$\bff_4$ is a relatively minimal rational surface with the
exceptional section $s$ such that $s^2=-4$. We denote
$E=\pi^{-1}(s)$ and by $e$ its class. The elliptic pencil $|C|$ is
the preimage by $\pi$ of the the rational pencil of $\bff_4$. We
denote its class by $c$. Then $c^2=0$, $e^2=-2$ and $c\cdot e=1$.
Since $mP\cdot c=2$ and $P\cdot e>0$, we get two cases:

{\bf Case $(\bff_4)^{(1)}$:} $m=2$ and $\wP=2P$ where
$P=n_1c+e$ where $n_1\ge 3$. Then
$n=P^2=2n_1-2$. We have $n\equiv 0\mod 2$ and $n\ge 6$.

{\bf Case $(\bff_4)^{(2)}$:} $m=1$ and $\wP=P=n_1C+2E$ where $n_1\ge 5$ and
$n_1\equiv 1\mod 2$. Then
$n=P^2=4n_1-8$. We have $n\equiv 4\mod 8$ and $n\ge 12$.

Using Riemann-Roch Theorem on $\bff_4$, it is easy to see that in both these
cases the linear system $|\wP|$ is the preimage by $\pi$ of the corresponding
linear system from $\bff_4$. It follows that $|\wP|:X\to Y\subset \bp^N$ is 
the quotient map.

\medskip

{\bf Case $\bff_0=\bp^1\times \bp^1$:} $S\cong U(2)$ and $Y=X/\{1,\tau\}\cong
\bp^1\times \bp^1$. We denote by $e_1$ and $e_2$ classes of preimages of
$\text{pt}\times \bp^1$ and $\bp^1\times \text{pt}$. We have
$e_1^2=e_2^2=0$ and $e_1\cdot e_2=2$. Then $C$ has the class
of $e_1$ or $e_2$, and we obtain two cases. Either
$\wP=P=n_1e_1+e_2$, $n_1\ge 1$ (case $(\bff_0)_{(1)}$) or
$\wP=P=e_1+n_1e_2$, $n_1\ge 1$ (case $(\bff_0)_{(2)}$).
Then $n=P^2=4n_1$ where $n\equiv 0\mod 4$
and $n\ge 4$. In both cases the linear system $|\wP|$
is the preimage of the corresponding linear system from $\bff_0$.
It follows that $|\wP|:X\to Y\subset \bp^N$ is the quotient map.
(Over $\br$, it
is important to distinguish the cases $(\bff_0)_{(1)}$ and $(\bff_0)_{(2)}$
which are symmetric over $\bc$.)

\medskip

{\bf Case $\bff_2$:} $S\cong U(2)$ and $Y=X/\{1,\tau\}\cong \bff_2$.
This case can be
considered as the degeneration of the previous one. We denote by $e$
the class of the preimage of the exceptional section of the relatively minimal
rational surface $\bff_2$ (it is union of two disjoint non-singular rational
curves on $X$ which are conjugate by $\tau$; then the Picard lattice of $X$
has the rank at least three),
and by $c$ the class of the preimage of the rational pencil of $\bff_2$.
We have $e^2=-4$, $c^2=0$ and $c\cdot e=2$. Then $C$ has the class $c$,
and $\wP=P=n_1c+e$ where $n_1\ge 3$. Then $n=4n_1-4$ where $n\equiv 0\mod 4$
and $n\ge 8$. Again $|\wP|$ is the preimage of the corresponding
linear system from $\bff_2$ and $|\wP|:X\to Y\subset \bp^N$
is the quotient map.

\medskip

These considerations show that classification of hyper-elliptically
polarized K3 surfaces $(X,\wP)$ is equivalent to classification of
K3 surfaces with non-symplectic involution $(X,\tau)$ when $Y=X/\{1,\tau\}$ is
$\bp^2$ or $\bff_r$ for $r=0,1,2,4$, by picking up $|\wP|=\pi^\ast |Q|$
where $|Q|$ is the appropriate linear system of $Y$. Considering their
deformations gives the corresponding polarized K3 surfaces we are looking for.
The following statement shows that we can drop $\bff_2$ from the consideration.

\begin{lemma}
\label{F2toF0}
Any real polarized K3 surface which is a deformation of
a real hyper-elliptically polarized K3 surface coming from $\bff_2$ is also
a deformation of a real hyper-elliptically polarized K3 surface coming from
$\bff_0$.
\end{lemma}

\Proof Let $(X,\wP)$ be a real hyper-elliptically polarized K3 surface where
$\wP\in S\cong U(2)$ comes from $Y=X/\{1,\tau\}$ where $Y$ is $\bff_2$. Let
$\varphi$ be the corresponding anti-holomorphic involution on $X$ which 
commutes 
with $\tau$. Thus, the triplet $(X,\tau,\varphi)$ is the
real K3 surface with the non-symplectic involution $\tau$. If $Y\cong \bff_2$,
the Picard lattice $S_X$ contains $S$, but it has the rank at least 3.
Changing complex structure of $X$ a little and applying Global Torelli Theorem
for K3 surfaces \cite{PS71}, we can construct another real K3 surface with a
non-symplectic involution 
$(\widetilde{X},\widetilde{\tau}, \widetilde{\varphi})$
which has the same action of involutions on $H_2(\widetilde{X},\bz)$,
but has the Picard lattice $S_{\widetilde{X}}=S$ of the rank two.
Then $\widetilde{X}/\{1,\widetilde{\tau}\}\cong \bff_0$,
and the same class $\wP$ comes from $\bff_0$. Both these real polarized K3
surfaces give the same integral polarized K3 involutions
$(L,\varphi=\widetilde{\varphi},\wP)$ and, by Theorem \ref{intpolinvtheorem},
the same connected component of moduli of real polarized K3 surfaces.
This finishes the proof.
\QED

\medskip

Hyper-elliptically
polarized K3 surfaces $(X,\wP)$ where $\wP$ comes from $Y=X/\{1,\tau\}$ where
$Y\cong \bff_r$, $r=0,1,4$, were called in \cite{NikulinSaito05} as 
{\it general 
K3 double rational scrolls} because $|\wP|$ defines an embedding
$\bff_r=Y\subset \bp^N$ as a rational scroll
(the degree of $Y$ is equal to $N-1$, see \cite{Saint-Donat74}),  and $X$
is its double covering having in general the smallest possible Picard number 2.
By above
considerations and Lemma \ref{F2toF0},
deformations of real hyper-elliptically polarized K3 surfaces (which we want to
classify) are equivalent to deformations of general real K3 double 
rational scrolls.

In \cite{NikulinSaito05}, classification of connected components of moduli
of general real K3 double rational scrolls had been obtained. Their K3
polarized deformations had been studying.
Using results of Sect. \ref{secmoduliK3}, here want to finalize 
the results about deformations.

\medskip

The main result of \cite{NikulinSaito05} about deformations is

\begin{theorem} (see Theorem 30 in \cite{NikulinSaito05}) A genus invariant
\eqref{genK3} of a real polarized K3 surface can be obtained as a
deformation of some general real K3 double rational scroll (equivalently, of
a real hyper-elliptically polarized K3 surface) if and only if the following
condition is valid:
$$
n\le 4\ if\ either\ (r,a)=(20,2)\ or\ r+a=22\ and\ \delta_{\varphi P}=0.
$$
Equivalently,
$$
n\le 4\ if \ either\ X(\br)=(T_0)^{10}\ or\ X(\br)=(T_0)^{k}\ and\
X(\br)\sim P\mod 2\ in\ H_2(X(\bc),\bz).
$$
\label{theoremgendeform}
\end{theorem}

\medskip

Since for $n=2$ or $n=4$, the genus invariants give 
only one connected component
(see Sect. \ref{secmoduliK3}), we then obtain from Theorem 
\ref{theoremdeformgen} 
the following result which is trivial for $n=2$ and  
had been proved for $n=4$ in
\cite{NikulinSaito05}.

\begin{theorem} (see \cite{NikulinSaito05})
\label{n=2n=4}
Any real polarized K3 surface $(X,P^\prime)$ of the primitive degree
$n=2$ or $n=4$ is a deformation of
a general real K3 double rational scroll (equivalently, of a real 
hyper-elliptically polarized K3 surface).
\end{theorem}

Thus, further we can assume that $n\ge 6$. Applying the connectedness
Theorem \ref{theoremuniqueness},
from Theorem \ref{theoremgendeform} we get the result.

\begin{theorem} Let $(X,P^\prime)$ be a real polarized K3 surface of
the primitive degree $n\ge 6$,  and $X(\br)\not= T_1\amalg (T_0)^8,\ (T_0)^9$.

Then $X$ is a deformation of a general real K3 double rational scroll
(equivalently of a real hyper-elliptically polarized K3 surface) if and only if
$X(\br)\not=(T_0)^{10}$ and ($X(\br)\not\sim P\mod 2$ in $H_2(X(\bc),\bz)$ when
$X(\br)=(T_0)^k$).
\label{theoremdeformgen}
\end{theorem}

\medskip

Now let us assume that $n\ge 6$ and $X(\br)=T_1\amalg (T_0)^8$. Equivalently,
$(r,a,\delta_\varphi)=(19,1,1)$. Then
we should apply invariants of Theorems \ref{theorem1911} and \ref{theorem1910}.

By results of \cite{NikulinSaito05} (especially see Sect. 8.1 in 
\cite{NikulinSaito05}), 
there are three types of the deformations, and
they can be described symbolically  as shown (see also considerations above).

{\it The Case of $(\bff_4)^{(1)}$:} $n\in (\bff_4)^{(1)}$ if and only if
$n\equiv 0\mod 2$ and $n\ge 4$; $P=(n/2+1)C+E$;
$$
D_n:((\bff_4)^{(1)};r=19,a=1,H=0,\delta_{\varphi S}=1)
\Longrightarrow
(n;r=19,a=1,\delta_P=1,\delta_\varphi=1, \delta_{\varphi P}=1).
$$

In this case
$$
L=[C,E]\oplus [g_1,g_2;(g_1+g_2)/2]\oplus U\oplus 2E_8
$$
where $C^2=0$, $E^2=-2$ and $C\cdot E=1$, $g_1^2=2$, $g_2^2=-2$ and
$g_1\perp g_2$. The involutions $\tau$ and $\varphi$ are characterized by
their eigenvalue 1 parts $L^\tau=[C,E]$ and 
$L^\varphi= [g_2]\oplus U\oplus 2E_8$. 
Then the triplet $(L,\tau,\varphi)$ has the required invariants.

The lattice $L^{\varphi,P}=\bz P\oplus \bz g_2\oplus U\oplus 2E_8$
is $\langle n\rangle\oplus
\langle -2 \rangle$ modulo the unimodular lattice $U\oplus 2E_8$. Here $P$ and
$g_2$ are the standard generators of $\langle n \rangle$ and 
$\langle -2 \rangle$ 
respectively. Its orthogonal complement in $L$ is
$L_{\varphi,P}=\bz Q\oplus \bz g_1$ where $Q=(n/2-1)C-E$ and $Q^2=-n$.
It follows that $L_{\varphi,P}=\langle -n\rangle\oplus \langle 2 \rangle$
where $Q$ and $g_2$ are the standard generators of $\langle -n\rangle$ and
$\langle 2\rangle$ respectively.

We have $Q/n+P/n=C\in L$ and $g_2/2+g_1/2=(g_1+g_2)/2\in L$. It follows
that the connected component of moduli is standard.

\bigskip

{\it The Case of $(\bff_4)^{(2)}$:} $n\in (\bff_4)^{(2)}$ if and only if
$n\equiv 4\mod 8$ and $n\ge 12$; \linebreak $P=(n/4+2)C+2E$;
$$
D_n:((\bff_4)^{(2)};r=19,a=1,H=0,\delta_{\varphi S}=1)
\Longrightarrow
(n;r=19,a=1,\delta_P=1,\delta_\varphi=1, \delta_{\varphi P}=1).
$$

In this case, everything is the same as in the previous one, only 
$P=(n/4+2)C+2E$ and $Q=(n/4-2)C-2E$. Thus,
$$
\frac{P}{n/2}+\frac{Q}{n/2}=C\in L,\ \
\frac{P}{4}+\left(-\frac{Q}{4}\right)=C+E\in L
$$
where $n\equiv 4\mod 8$.
It follows that the connected component of moduli is different from the
standard one by $(-1)$ in the $2$-component of the discriminant form.

\bigskip

{\it The Case of $\bff_1$:} $n\in \bff_1$ if and only if
$n\equiv 2\mod 4$ and $n\ge 6$; $P=\left((n+2)/4\right)h+
\left((2-n)/4\right)e$; 
$$
D_n:(\bff_1;r=19,a=1,H=[h],\delta_{\varphi S}=0,v=h)
\Longrightarrow
(n;r=19,a=1,\delta_P,\delta_\varphi=1, \delta_{\varphi P})\ \text{where}
$$
$$
\delta_P=\delta_{\varphi P}=
\left\{\begin{array}{cl}
0 &\mbox{if}\ n\equiv 2\mod 8\\
1 &\mbox{if}\ n\equiv -2\mod 8
\end{array}\right .  \ .
$$

In this case, the lattice
$$
L=[h,e,g_1,g_2;(h+g_2)/2,(e+g_1)/2]\oplus U\oplus 2E_8
$$
where $h,e,g_1,g_2$ are orthogonal to each other and $h^2=2$,
$e^2=-2$, $g_1^2=2$, $g_2^2=-2$.
The involutions $\tau$ and $\varphi$ are characterized by their 
eigenvalue 1 parts $L^\tau=[h,e]$ and 
$L^\varphi$=$\bz g_2\oplus U\oplus  2E_8$.
Then the triplet $(L,\tau, \varphi)$ has the required invariants.

Assume that $n\equiv -2\mod 8$. Then
$L^{\varphi,P}=\bz P \oplus \bz g_2 \oplus U\oplus 2E_8$ which is
$\langle n\rangle\oplus
\langle -2 \rangle$ modulo the unimodular lattice $U\oplus 2E_8$.
Here $P$ and $g_2$ are the standard generators of 
$\langle n \rangle$ and $\langle -2 \rangle$ 
respectively.
Its orthogonal complement in $L$ is $L_{\varphi,P}=\bz Q \oplus \bz g_1$
where $Q=((n-2)/4)h-((n+2)/4)e$ and $Q^2=-n$. Thus,
$L_{\varphi,P}= \langle -n \rangle \oplus \langle 2 \rangle$ where
$P$ and $g_1$ are the standard generators of $\langle -n \rangle$ and
$\langle 2 \rangle$ respectively.

We have
$$
e=\frac{n-2}{2n}P-\frac{n+2}{2n}Q,\ \ h=\frac{n+2}{2n}P+\frac{2-n}{2n}Q.
$$
It follows
$$
\frac{P}{n/2}+\frac{Q}{n/2}=h-e\in L,\ \
$$
and
$$
-\frac{n}{2}\left( \frac{g_2}{2} \right)+\frac{n-2}{4}\left(\frac{Q}{2}\right)
=-\frac{n}{2} \left( \frac{h+g_2}{2} \right) + \frac{n+2}{8}P\in L.
$$
This shows that the connected component of moduli is different from the
standard one by the non-trivial automorphism of the $2$-component of the
discriminant form which has the automorphism group $\bz/2$ in this case.

If $n\equiv 2\mod 8$, then the $2$-component of the discriminant form is
trivial. In this case, our (or the same) calculations show that
the connected component is standard.

\medskip

Thus, we finally get the result.

\begin{theorem} Let $(X,P^\prime)$ by a real polarized K3 surface of
the primitive degree $n\ge 6$, and $X(\br)=T_1\amalg (T_0)^8$.

Then $X$ is a deformation of a general real K3 double rational scroll 
(equivalently, 
of a real hyper-elliptically polarized K3 surface) if and only if
the following conditions satisfy

(a) if $n\equiv 0,\ 2\mod 8$, then $(X,P^\prime)$ is standard;

(b) if $n\equiv 4,\ 6\mod 8$, then $(X,P^\prime)$ is different from a standard
only over $2$.
\label{theoremT1T09}
\end{theorem}

\bigskip

Now let us assume that $n\ge 6$ and $X(\br)=(T_0)^9$. Equivalently,
$(r,a,\delta_\varphi)=(19,3,1)$. Then
we should apply invariants of Theorems \ref{theorem1931} since
$\delta_{\varphi P}=1$.

By results of \cite{NikulinSaito05}
(especially see Sect. 8.1 in \cite{NikulinSaito05}),
there are two types of the deformations, and
they can be described symbolically  as shown (see also considerations above).

{\it The Case of hyperboloid $\bhh_{(1)}$:} $n\in \bhh_{(1)}$
if and only if $n\equiv 0\mod 4$ and $n\ge 4$; $P=(n/4)e_1+e_2$;
$$
D_n:(\bhh_{(1)};r=19,a=3,H=[e_1,e_2],\delta_{\varphi S}=1)
\Longrightarrow
(n;r=19,a=3,\delta_P=0,\delta_\varphi=1,
\delta_{\varphi P}=1).
$$

In this case
$$
L=[e_1,e_2,e_1^\prime,e_2^\prime;(e_1+e_1^\prime)/2,(e_2+e_2^\prime)/2]
\oplus [g, g^\prime;(g+g^\prime)/2]\oplus 2E_8
$$
where $e_1^2=e_2^2=0$, $e_1\cdot e_2=2$; $(e_1^\prime)^2=(e_2^\prime)^2=0$,
$e_1^\prime\cdot e_2^\prime=-2$; $[e_1,e_2]\perp[e_1^\prime,e_2^\prime]$;
$g^2=2$, $(g^\prime)^2=-2$, $g\perp g^\prime$.
The involutions $\tau$ and $\varphi$ are characterized by 
$L^\tau=[e_1,e_2]$ and
$L^\varphi=[e_1^\prime,e_2^\prime]\oplus [g^\prime]\oplus 2E_8$.
Then the triplet $(L,\tau,\varphi)$ has the required invariants.

The lattice $L_\varphi=[e_1,e_2,g]$. We have $\bz P=\langle n \rangle$
where $P=(n/4)e_1+e_2$ is the standard generator of $\langle n \rangle$.
Its orthogonal complement in
$L_\varphi$ is $L_{\varphi,P}=\bz Q\oplus \bz g$ where 
$Q=(n/4)e_1-e_2$ and $Q^2=-n$. 
Thus $L_{\varphi,P}=\langle -n \rangle\oplus \langle 2 \rangle$ where
$Q$ and $g$ are the standard generators of $\langle -n\rangle$ and 
$\langle 2 \rangle$
respectively. We have
$$
\frac{Q}{n/2}+\frac{P}{n/2}=e_1\in L.
$$
It follows (see Theorem \ref{theorem1931}) that the connected component 
of moduli is standard.

\medskip

{\it The Case of $\bff_1$:} $n\in \bff_1$ if and only if
$n\equiv 2\mod 4$ and $n\ge 6$; $P=\left((n+2)/4\right)h+
\left((2-n)/4\right)e$; 
$$
D_n:(\bff_1;r=19,a=3,H=[h,e],\delta_{\varphi S}=1)
\Longrightarrow
(n;r=19,a=3,\delta_P=0,\delta_\varphi=1, \delta_{\varphi P}=1).
$$

In this case
$$
L=[h,e,h^\prime,e^\prime;(h+h^\prime)/2,(e+e^\prime)/2]
\oplus [g, g^\prime;(g+g^\prime)/2]\oplus 2E_8
$$
where $h^2=2$, $e^2=-2$, $h\perp e$; $(h^\prime)^2=-2$, $(e^\prime)^2=2$,
$h^\prime\perp e^\prime$; $[h,e]\perp [h^\prime,e^\prime]$;
$g^2=2$, $(g^\prime)^2=-2$, $g\perp g^\prime$.
The involutions $\tau$ and $\varphi$ are characterized by $L^\tau=[h,e]$ and
$L^\varphi=[h^\prime,e^\prime]\oplus [g^\prime]\oplus 2E_8$.
Then the triplet $(L,\tau,\varphi)$ has the required invariants.

The lattice $L_\varphi=[h,e,g]$. We have $\bz P=\langle n \rangle$
where $P=\left((n+2)/4\right)h+\left((2-n)/4\right)e$ is the standard
generator of $\langle n \rangle$. Its orthogonal complement in
$L_\varphi$ is $L_{\varphi,P}=\bz Q\oplus \bz g$ where
$Q=((n-2)/4)h-((n+2)/4)e$ and $Q^2=-n$.
Thus $L_{\varphi,P}=\langle -n \rangle\oplus \langle 2 \rangle$ where
$Q$ and $g$ are the standard generators of 
$\langle -n\rangle$ and $\langle 2 \rangle$ 
respectively. We have
$$
\frac{Q}{n/2}+\frac{P}{n/2}=h-e\in L.
$$
It follows (see Theorem \ref{theorem1931}) that the connected 
component of moduli is standard.

Thus, we obtain the result.

\begin{theorem} Let $(X,P^\prime)$ by a real polarized K3 surface of
a primitive degree $n\ge 6$, and $X(\br)=(T_0)^9$.

Then $X$ is a deformation of a general real K3 double rational scroll (i. e.
of a real hyper-elliptically polarized K3 surface) if and only if 
$X(\br)\not\sim P$ 
in $H_2(X(\bc),\bz)$ and $(X,P^\prime)$ is standard.
\label{theoremT09}
\end{theorem}

Let us unify all these results in one final statement.

\begin{theorem}
\label{theoremhypK3}
A real polarized K3 surface $(X,P^\prime)$ is
a deformation of a general real K3 double rational scroll 
(equivalently, of a real 
hyper-elliptically polarized K3 surface) if and only if one 
of conditions (i)---(iv) 
below satisfies:

(i) The primitive degree $n=2$ or $4$.

(ii) The primitive degree $n\ge 6$, and
$X(\br)\not=T_1\amalg (T_0)^8,\ (T_0)^9,\ (T_0)^{10}$, and
$X(\br)\not\sim P\mod 2$ in $H_2(X(\bc),\bz)$ if $X(\br)=(T_0)^k$.

(iii) The primitive degree $n\ge 6$, and
$X(\br)=T_1\amalg (T_0)^8$, and $(X,P^\prime)$ is standard if
$n\equiv 0,\ 2\mod 8$, and $(X,P^\prime)$ is different from standard 
only over $2$ 
if $n\equiv 4,\ 6\mod 8$.

(iv) The primitive degree $n\ge 6$, and
$X(\br)=(T_0)^9$, and $X(\br)\not\sim P\mod 2$ in $H_2(X(\bc),\bz)$, and
$(X,P^\prime)$ is standard.
\end{theorem}


\end{document}